\documentclass[12pt,reqno]{amsart}
\usepackage{amsmath}
\usepackage{amssymb}
\newtheorem{defn0}{Definition}[section]
\newtheorem{prop0}[defn0]{Proposition}
\newtheorem{thm0}[defn0]{Theorem}
\newtheorem{lemma0}[defn0]{Lemma}
\newtheorem{corollary0}[defn0]{Corollary}
\newtheorem{example0}[defn0]{Example}
\newtheorem{remark0}[defn0]{Remark}
\newtheorem{conjecture0}[defn0]{Conjecture}

\newenvironment{proposition}{\bigskip \begin{prop0}}{\end{prop0}}
\newenvironment{theorem}{\bigskip \begin{thm0}}{\end{thm0}}
\newenvironment{lemma}{\bigskip \begin{lemma0}}{\end{lemma0}}

\newenvironment{remark}{ \begin{remark0}\rm}{\end{remark0}}

\newcommand{\propref}[1]{Proposition~\ref{#1}}

\newcommand{\remref}[1]{Remark~\ref{#1}}

\def\mmax{{\bf m}}                   
\def\res{{\bf k}}                   
\def\oI{{o_{\max}(I)}}
\def\length{{ \mathrm{Length}_R}}

\begin{document}

\title[
Upper bounds of Hilbert coefficients and Hilbert functions ] {{\bf
Upper bounds of Hilbert coefficients and Hilbert
functions }}

\author[Juan Elias]{Juan Elias ${}^{*}$}
\thanks{
${}^{*}$Partially supported by DGICYT BFM2001-3584}
\address{Departament d'\`{A}lgebra i Geometria
\newline \indent Facultat de Matem\`{a}tiques
\newline \indent Universitat de Barcelona
\newline \indent Gran Via 585, 08007
Barcelona, Spain} \email{{\tt elias@ub.edu}}

\begin{abstract}
Let $(R,\mmax)$ be a $d$-dimensional Cohen-Macaulay local ring.
 In this note we prove, in a
very elementary way, an upper bound of the first normalized Hilbert
coefficient of a $\mmax$-primary ideal $I\subset R$ that improves
all known upper bounds unless for a finite number of cases, see
\remref{unless}. We also provide new upper bounds of the Hilbert
functions of $I$ extending the known bounds for the maximal ideal.
\end{abstract}

\maketitle

\baselineskip 14.5pt

\section*{Introduction}

The aim  of this paper is to provide new upper bounds for the first
normalized  Hilbert coefficient and the Hilbert function of a
$\mmax-$primary ideal. In the first  section  we give an upper bound
of the first normalized Hilbert coefficient $e_1(I)$ for an
$\mmax$-primary ideal $I$ taking account the multiplicity of $I$ and
the Hilbert polynomial of $R$, \propref{newb}. This new bound
improves the bound of \cite{RV05} eventually unless a finite number
of cases. For instance: if $I\subset \mmax^2$ then the bound of this
paper is better than the bound of \cite{RV05}, see  \remref{unless}.
The advantage of the result proved here is that its proof is very
elementary skipping the theory of Ratliff-Rush closure or the
integral closure of ideals. For more details on the upper bounds of
the first Hilbert coefficient see the beginning of section 1.

In the second section we give several upper bounds of the Hilbert
function of $I$ generalizing some results of \cite{RVV01}. The main
ingredients of these bounds are the upper bound for the first
normalized Hilbert coefficient and the upper bound of the
postulation number, \propref{sta}. This upper bound of the
postulation is given
 in terms of the normalized Hilbert coefficients of $I$ and
the length of the conormal module $I/I^2$.

\medskip
\noindent {\sc Notations.} Let $R$ be a $d$-dimensional
Cohen-Macaulay ring with maximal ideal $\mmax$. Without loss of
generality we can assume that the residue field $\res=R/\mmax$ is
infinite.
 Let $I$ be an $\mmax$-primary ideal of $R$,
the Hilbert-Samuel  function of $I$ is defined by
$H^{0}_{I}(n)=\length(I^n/I^{n+1})$, $n\ge 0$. The $i$-th
Hilbert-Samuel  function of $I$, $i\ge 1$, is
$H^{i}_{I}(n)=\sum_{j=0}^n H^{i-1}_{I}(j)$, $n\ge 0$. It is well known that
there exist integers $e_j(I) \in \mathbb Z$ such that
$$
h^i_{I}(X)=\sum_{j=0}^{d+i-1}
(-1)^j e_j(I) \binom{X+d+i-j-1}{d+i-j-1}
$$

\noindent is the $i$-Hilbert-Samuel polynomial of $I$, i.e.
$H^i_{I}(n)=h^i_{I}(n)$ for  $n \gg 0$, $i\ge 0$.
 We denote by
$pn(I)$ the postulation number of $I$, i.e. the least integer
$pn(I)$ such that $H^0_{I}(n)=h^0_{I}(n)$ for all $n \ge pn(I)$; the
stability $st(R)$ of $R$  is the supreme of $pn(I)$ where $I$ range
over all the $\mmax-$primary ideals of $R$. The order or initial
degree of $I$ is the integer $o_{\mmax}(I)$ such that $I\subset
\mmax^{\oI}$ and $I\nsubseteq \mmax^{\oI+1}$, i.e. the last integer
$n$ such that $H^0_I(n)=H^0_R(n)$. We denote by
$b(I)=H^0_I(1)=\length(I/I^2)$ the length of the conormal module of
$I$. In the maximal ideal case we write $H^*_R=H^*_{\mmax}$,
$b(R)=b(\mmax)$ is the  embedding dimension of $R$, and
$e_i(R)=e_i(\mmax)$, $i\ge 0$.

The $i$-Poincar{\'e} series of $I$ is the generating function of
$H^i_I$:
$$
P_I^i= \sum_{n\ge 0} H^i_I(n) \; T^n \; \; \in {\mathbb Z}[[T]]
$$
$i \ge 0$. Since $H^i_I$ is asymptotically polynomial there exists
a degree $pn(I)+d-1$ polynomial $f(T)\in \mathbb Z [T]$ such that
$$
P_I^i= \frac{f(T)}{(1-T)^{d+i}}
$$

An ideal $J\subset R$ is a reduction of $I$ if
there exists an integer $r$ such that $I^{r+1}=J I^r$. The
reduction number of $I$ is the least integer $r(I)$ for which
there exists a reduction $J$ of $I$ such that $I^{r(I)+1}=J I^{r(I)}$.
 We
denote by $\mu(I)$ the minimal number of generators of $I$.

\medskip
The author acknowledges the useful comments of M.E. Rossi.

\medskip
\section{Upper bound of $e_1(I)$}

In \cite{Eli90} we characterize the Hilbert-Samuel  polynomials of
one-dimensional Cohen-Macaulay local rings $R$,  in particular we
give a sharp upper bound of the first normalized Hilbert coefficient
$e_1(R)$ in terms of the multiplicity of $R$ and its embedding
dimension $b(R)$, see also \cite{Eli01}. This upper bound was
naturally generalized in \cite{Eli05} to $\mmax$-primary ideals $I$.
In  \cite{RV05} Rossi and Valla improved the bound of \cite{Eli05}
by considering the minimal number of generators of $I$ and the
length of $R/I$.

In this section  we give
 an  upper bound of the first normalized Hilbert coefficient
of $I$ in terms of the Hilbert-Samuel polynomial of $R$ and the
multiplicity of $I$ improving the bound of Rossi and Valla.
 We include  some known results of Northcott
completing the picture on $e_1(I)$.

\begin{proposition}
\label{newb} Let $R$ be  a $d$-dimensional Cohen-Macaulay ring. Let
$I$ be a $\mmax$-primary ideal of $R$, then
\begin{enumerate}
\item[(i)] {\rm [Northcott]} $0\le e_0(I)- \length(R/I)\le
e_1(I)$, \item[(ii)] $e_0(R) \oI \le e_0(I)$, \item[(iii)] If
$d=1$ then $e_1(I)\le pn(I)(e_0(I)- e_0(R) \oI ) +
e_1(R).$\newline For $d\ge 1$ it holds
$$e_1(I)\le (e_0(R)-1) (e_0(I)- e_0(R) \oI ) + e_1(R).$$
\end{enumerate}
\end{proposition}
\begin{proof}
$(i)$ is due to Northcott, \cite{Nor60}. $(ii)$ and $(iii)$ for
$d=1$. Since the residue field is infinite we can assume that $R$ is
one-dimensional. Notice that from \cite{Lip71}, Theorem 1.5, and the
fact $I \subset \mmax^{\oI}$ we deduce
\begin{multline*}
e_0(R)(n\, \oI)- e_1(R) \le \length(R/\mmax^{n\, \oI})\le\\
\le \length(R/I^n)= e_0(I) n -e_1(I)
\end{multline*}
for all $n\ge pn(I)$. Hence  we get
$
e_0(R) \oI \le e_0(I),
$
and for $n = pn(I)$ we deduce
$$
e_1(I)\le  pn(I) (e_0(I)- e_0(R) \oI ) + e_1(R).
$$

\noindent $(iii)$ for  $d\ge 1$. Notice that, since the residue field is infinite, $e_1(I)$ and
the integer
$$
(e_0(R)-1) (e_0(I)- e_0(R) \oI ) + e_1(R)
$$
is stable by quotienting $R$ by a generic  regular sequence. Hence
we may assume that $R$ is one-dimensional. Notice that $pn(I)\le
e_0(R)-1$, \cite{SV74} Theorem 2.5, so from the
inequality $(iii)$ in the one-dimensional case we get
$$
e_1(I)\le (e_0(R)-1)  (e_0(I)- e_0(R) \oI ) + e_1(R).
$$
\end{proof}

\medskip
Let $R$ be a $d$-dimensional Cohen-Macaulay local ring and let $I$
be an $\mmax$-primary ideal of $R$. We denote by  $\epsilon(I)$ the
upper bound of $e_1(I)$ of the last result:
$$\epsilon(I)=(e_0(R)-1)(e_0(I)- e_0(R) \oI ) + e_1(R).$$

\noindent
Rossi and Valla gave the following  upper bound of $e_1(I)$, \cite{RV05},
$$
e_1(I)\le \rho(I)=\binom{e_0(I)}{2}- \binom{\mu(I)-d}{2}-\length(R/I)+1.
$$

\medskip
In the next result we prove that the bound $\epsilon(I)$  improves
 $\rho(I)$ unless a finite number of cases.

\begin{proposition}
Let $R$ be  a $d$-dimensional Cohen-Macaulay ring. Let $I$ be a
$\mmax$-primary ideal of $R$ such that $e_0(I)\ge 2 e_0(R)-1$ then
$$
\epsilon(I) \le  \rho(I).
$$
If $b(R)\ge 3$ or $\oI\ge 3$ then the last inequality is strict.
\end{proposition}
\begin{proof}
Since the residue field is infinite we can consider that $d=1$. A
simple computation and by using the inequalities $\length(R/I)\le
e_0(I)$ and
$$
e_1(R)\le \binom{e_0(R)}{2} - \binom{b-1}{2},
$$
$b=b(R)$, \cite{Eli90}, yields us
\begin{multline*}
2(\rho(I)- \epsilon(I))\ge       \qquad  \qquad  \qquad           \\
\qquad \qquad \ge e_0(I)(e_0(I) - 2\, e_0(R)-1)+ 2\, e_0(R)(e_0(R)(\oI -1)- \oI +2)+ 2 \binom{b-1}{2}.
\end{multline*}

\noindent
Let us assume $\oI =1$. Then we have
$$
2(\rho(I)- \epsilon(I))\ge  e_0(I)(e_0(I) - 2\, e_0(R)-1)+ 2\,
e_0(R) + 2 \binom{b-1}{2}.
$$
If $e_0(I)\ge 2 \, e_0(R)-1$ then we deduce
$$
2(\rho(I)- \epsilon(I))\ge  \, e_0(R)(-2)+ 2\, e_0(R)+ 2
\binom{b-1}{2}= 2 \binom{b-1}{2} \ge 0.
$$
If $b(R)\ge 3$ then $2(\rho(I)- \epsilon(I))> 0$.

\noindent
Let us assume that $\oI \ge  2$.
Since $e_0(I) \ge \oI \, e_0(R)$ we have
$$
2(\rho(I)- \epsilon(I)) \ge   e_0(R)((\oI^2 - 2) e_0(R) - 3 \, \oI
+ 4)+ 2 \binom{b-1}{2}.
$$
Notice that $\oI^2 - 2 \ge 0$ and $e_0(R)\ge 1$, so
$$
2(\rho(I)- \epsilon(I)) \ge   \oI^2 - 3 \, \oI + 2 + 2
\binom{b-1}{2}\ge 0.
$$
If $b(R)\ge 3$ or $\oI\ge 3$ then $2(\rho(I)- \epsilon(I))>0$.
\end{proof}

\begin{remark}
\label{unless} Notice that the non covered cases by the last
proposition are finite: $\oI = 1$ and   $e_0(R) \le e_0(I) \le 2
e_0(R)-2.$ In particular if $I\subset \mmax^2$ then $\epsilon(I)\le
\rho(I)$.
\end{remark}

\begin{remark} {\bf  Comparison with other  upper bounds.}
Let us assume that $R$ is one dimensional. In the proof of Theorem
3.1 of \cite{RVV01} it is proved that
$$
e_1(I)\le \alpha(I)=e_0(I)-\length(R/I)+
pn(I)(e_0(I)-2)-\binom{pn(I)+1}{2}
$$
provided that $e_0(I)\neq e_0(R)$. It is easy to see that this
bound is not comparable with the bound of \propref{newb} $(iii)$
$$
\beta(I)= pn(I)(e_0(I)- e_0(R) \oI ) + e_1(R),
$$
i.e. there exists ideals $I$ such that $e_1(I)\le \alpha(I)<
\beta(I)$, resp. $e_1(I)\le\beta(I)< \alpha(I)$.

Since $pn(I)$ behaves bad under quotienting then by considering the
inequality $pn(I)  \le e_o(I)-1$ and the upper bound $\alpha(I)$,
both in the one-dimensional case, Rossi and Valla prove
for any $d$-dimensional local ring
$R$, $d \ge 1$, that
$$
e_1(I)\le \binom{e_0(I)-2}{2}
$$
provided $e_0(I)\neq e_0(R)$,  \cite{RVV01} Proposition 3.1. An easy computation
shows that our bound is better:
$$
\epsilon(I)\le \binom{e_0(I)-2}{2}
$$
except a finite number of cases: $\oI=1$, $e_0(I)=e_0(R)+1$, and
$e_0(R)\le 2$ if $b(R)=e_0(R)+1$, or $e_0(R)\le 4$ if
$b(R)=e_0(R)$.

Recall that in \cite{Eli05}, Proposition 2.5, we prove that
$\epsilon(I)\le \binom{e_0(I)-s(I)}{2}$
if $e_0(I)$ is a prime number.
\end{remark}


\medskip
\section{Upper   bounds of the Hilbert function}

In this section we prove  several upper bounds of the Hilbert
function of $\mmax$-primary ideals  generalizing and improving the
upper bounds of the Hilbert function of $R$ given in  \cite{RVV01},
Theorem 3.4.
 The main
ingredients of the new  bounds are the upper  bound of  the first
normalized Hilbert coefficient  and a new upper bound
of the postulation number, \propref{sta}.

\medskip
Let $R$ be a one-dimensional Cohen-Macaulay
  local ring and let $I$ be an $\mmax-$primary ideal.
  For all $n\ge 0$ we set $v_n=e_0(I)-H^0_I(n)$.
  The set of integers $\{v_{\cdot}\}$ satisfies the following conditions:

\begin{enumerate}
  \item[(1)] $v_n \ge 0$, for all $n\ge 0$,
  \item[(2)]  $v_n=0$ if and only if $n\ge pn(I)$,
  \item[(3)]   $e_1(I)=\sum_{n\ge 0} v_n$.
\end{enumerate}
\noindent It is easy to prove from the definition  of the integers
$\{v_{\cdot}\}$ that, $n \ge 0$,
$$
H^1_I(n)=(n+1)e_0- \sum_{j=0}^n v_j.
$$

\noindent
 We define the integer $\lambda(I)=\length(R/I)$.
From the definition of $v_0$ and $v_1$ it is easy to see that
$v_0+v_1= 2 \, e_0(I)- \lambda(I) - b(I).$
Since $H^0_I(n)\le e_0(I)$  we can consider the integer
$$
\beta(I)=\left\{%
\begin{array}{ll}
    0 & \hbox{if }\,  e_0(I)=\lambda(I),\\
    1& \hbox{if } \, e_0(I)=b(I),\\
    e_1(I) - 2 e_0(I) + \lambda(I) + b(I)+2 & \hbox{if }\, e_0(I) >b(I).\\
\end{array}%
\right.
$$

\medskip
In the next result we first recall, without proof and for the readers convenience, a known
result on the stability due to Sally and Vasconcelos,   see \cite{SV74},
Theorem 2.5.
In the second part of the next result  we prove that $\beta(I)$ is also an upper bound of
the postulation number of $I$.

\begin{proposition}
\label{sta} Let $R$ be  a one-dimensional Cohen-Macaulay ring and
$I\subset R$ an  $\mmax-$primary ideal. Then
\begin{enumerate}
\item[(i)] {\rm [Sally-Vasconcelos]} $n(I)=pn(I) \le st(R)\le e_0(R)-1$,
\item[(ii)] $pn(I)\le \beta(I)$.
\end{enumerate}
\end{proposition}
\begin{proof} $(ii)$
Let us recall that if $H_I^0(n)=e_0(I)$ then $H_I^0(n+t)=e_0(I)$
for all $t\ge 0$, \cite{Lip71}, Lemma 1.8. From this result we get
the that $pn(I)\le \beta(I)$ in the cases $e_0(I)=\lambda(I)$ and
$e_0(I)=b(I)$.

Let us assume $e_0(I) >b(I)$, notice that we can write
$$
e_1(I)=\sum_{n\ge 0} v_n= v_0 + v_1 +\sum_{n\ge 2} v_n= 2\, e_0(I)
-\lambda(I) - b(I) + \sum_{n\ge 2} v_n,
$$
so we get
$$
\sum_{n\ge 2} v_n=e_1(I)- 2 \, e_0(I) + \lambda(I) + b(I) =
\beta(I)-2.
$$
From this inequality we deduce that $v_j=0$ for all $j\ge
\beta(I)$, in particular we get $pn(I)\le \beta(I)$.
\end{proof}

\begin{remark}
Rossi in  \cite{Ros00b}, Theorem 1.3, gives an upper bound of the
postulation number of an $\mmax$-primary ideal by a subtle
application of the determinantal trick:
$$
pn(I)\le w(I)=e_1(I)-e_0(I)+\lambda(I)+1,
$$
where $I$ is an $\mmax-$primary ideal of a $d-$dimensional
Cohen-Macaulay local ring $R$ with $d=1,2$. Let us remark that the
result of Rossi has its deep meaning in the two dimensional case.
Our bound, whose proof follows the steps of \cite{Ros00b}, improves
Rossi's bound, i.e. $\beta(I) \le w(I)$, for all $\mmax$-primary
ideal $I$ of a one-dimensional Cohen-Macaulay local ring $R$.
\end{remark}

\medskip
In the next results  we give  upper bounds
 of  the Hilbert function of $\mmax$-primary ideals, these results
generalize  \cite{RVV01}, Theorem 3.4.

\begin{theorem}
\label{RVV} Let $R$ be a one-dimensional Cohen-Macaulay local ring.
Let  $I$ be  an $\mmax-$primary ideal, then
$$
P^1_I =
\left\{%
\begin{array}{ll}
    \frac{\lambda(I)}{(1-T)^2}& \hbox{if }  \beta(I)=0\\ \\
    \frac{\lambda(I) + (b(I)-\lambda(I)) T
}{(1-T)^2} & \hbox{if  } \beta(I)=1\\
\end{array}%
\right.
$$
If $\beta(I)\ge 2$ then
$$
P^1_I \le
\frac{\lambda(I) + (b(I)-\lambda(I)) T + (e_0(I)-b(I)  -1)T^2
+ T^{\beta(I)}}{(1-T)^2}$$
\end{theorem}
\begin{proof}
We write $e_i=e_i(I)$, $\lambda=\lambda(I)$,  $\alpha=\alpha(I)$ and $\beta=\beta(I)$.

If $\beta=0$ then $H_I^0(n)=e_0$ for all $n\ge 0$,
in this case  we have $P^1_I = \lambda /(1-T)^2$.
If $\beta=1$ then $H_I^0(0)=\lambda$ and $H_I^0(n)=e_0$ for all $n\ge 0$, so
 $P^1_I = \lambda + (b-\lambda) T / (1-T)^2$.
We can assume that $\beta\ge 2$.
 First we compute the right hand side power series  of  the
inequality:

\begin{multline*}
\frac{\lambda + (b - \lambda) T + (e_0-b-1)T^2
+ T^{\beta}}{(1-T)^2}=\\
=\lambda +
\sum_{j=2}^{\beta-1}((j-1)(e_0-1)+b+\lambda) T^j+
\sum_{\beta}^{\infty}(e_0(j+1)-e_1) T^j.
\end{multline*}

\noindent
Since $H^1_I(0)=\lambda$ and $H^1_I(n) = e_0(n+1)-e_1$
for all $n \ge \beta$, \propref{sta}, we have to prove
$$
H^1_I(n) \le (n-1)(e_0-1)+ b +\lambda
$$
\noindent
for $n=2,\cdots, \beta -1$.
Notice that this inequality is equivalent to
\begin{equation*}
\tag*{(*)}
\sum_{j=0}^n v_j \ge 2 e_0 - b  - \lambda +n -1.
\end{equation*}

\noindent
If $v_n=0$ then $\sum_{j=0}^n v_j=e_1$ and we have to prove
$e_1 \ge  2 e_0 -  b - \lambda+n -1,$
this inequality is equivalent to  the assumption $n \le \beta -1$.

\noindent
If $v_n \neq 0$ then we may consider
$$
   \sum_{j=0}^n v_j  -n   = 2 e_0 - \lambda - b + \sum_{j=2}^n v_j  -n
      \ge 2 e_0 - \lambda - b -1
$$
\noindent
and we get (*).
\end{proof}

\medskip
We can provide the result for a $d$-dimensional local ring. For this end we have
to control the behavior of $b(I)$ by quotienting by a superficial sequence.
 Let $x_1,\cdots, x_{d-1} \in I$ be a superficial sequence of
$I$.
If we denote by $K=I/(x_1,\cdots,
x_{d-1})$ and $R^*=R/(x_1,\cdots, x_{d-1})$ we have
$\lambda(I)=\lambda(K)$,
$e_i(K)=e_i(I)$, $i=0,1$, and
$$
b(K)=b(I)-(d-1) \lambda(I) + \length(I^2\cap (x_1,\cdots, x_{d-1})/ (x_1,\cdots, x_{d-1}) I).
$$

\begin{theorem}
\label{RVVbis} Let $R$ be a $d$-dimensional Cohen-Macaulay local
ring. Let  $I$ be  an $\mmax-$primary ideal, and let $x_1,\cdots,
x_{d-1} \in I$ be a superficial sequence of $I$. If we denote by
$K=I/(x_1,\cdots, x_{d-1})$ then
$$
P^1_I \le
\left\{%
\begin{array}{ll}
    \frac{\lambda(K)}{(1-T)^2}& \hbox{if }  \beta(K)=0\\ \\
    \frac{\lambda(K) + (b(K)-\lambda(K)) T}{(1-T)^2} & \hbox{if  } \beta(K)=1\\ \\
\frac{\lambda(K) + (b(K)-\lambda(K)) T + (e_0(I)-b(K)  -1)T^2
+ T^{\beta(K)}}{(1-T)^2}  & \hbox{if } \beta(K)\ge 2\\
\end{array}%
\right.
$$
If the equality holds then $depth(gr_I(R)) \ge d-1$.
\end{theorem}
\begin{proof}
From Singh's inequality, \cite{Sin74},  we get
 $$
 P_I \le
\frac{P_{K}}{(1-T)^{d-1}},
 $$
\noindent
from \propref{RVV} we get the inequality of the claim.
If we have
$P_I=P_{K}/(1-T)^{d-1}$,
then $depth(gr_I(R)) \ge d-1$, \cite{Sin74}
\end{proof}

\medskip
We can skip the ideal $K$ of  the last proposition by considering a
new invariant attached to $I$ instead $\beta(I)$. Following
\cite{RVV01} we can  generalize to the $\mmax$-primary case some
results of that paper.

Let us recall that Rossi and Valla prove that
$$
e_1(I)\le \rho(I)= \binom{e_0(I)}{2}- \binom{\mu(I)-1}{2}-\length(R/I)+1
$$
According to this bound we define $\alpha(I)$ as the last integer
$n$ such that
$$\binom{n}{2} \le \binom{e_0(I)}{2} - e_1(I)-\length(R/I)+1,
$$
notice that  $\alpha(I) \ge \mu(I)-1$.
We define the integer
$$
\delta(I) = e_1(I) - 2 e_0(I) + (\alpha(I)+2) \; \lambda(I).
$$

\begin{lemma}
\label{millor} For all $\mmax$-primary ideal we have
$$
\beta(I)\le \delta(I)+2.
$$
\end{lemma}
\begin{proof}
We have to prove that
$$
e_1(I) - 2 e_0(I) +  \lambda(I) + b(I) +2=\beta(I)
\le \delta(I)+2=e_1(I) - 2 e_0(I) + (\alpha(I)+2) \; \lambda(I)+2,
$$
this inequality is equivalent to $ b(I) \le (\alpha(I)+1) \; \lambda(I).$
This inequality follows from the natural epimorphism of $R/I$-modules
$(R/I)^{\mu(I)}\rightarrow I/I^2$ and the fact $\alpha(I) \ge \mu(I)-1$.
\end{proof}

\medskip
The following result is a generalization to $\mmax$-primary ideals
of \cite{RVV01}, Theorem 3.4. In the one dimensional case the upper
bound of \propref{RVV} is termwise smaller that the upper bound in
$(i)$ of the next result.

\begin{proposition}
\label{RVV2} Let $R$ be a Cohen-Macaulay local ring of dimension
$d$. Let  $I$ be  an $\mmax-$primary ideal.
\begin{enumerate}
\item[(i)]  If $d=1$ then it holds
$$
P^1_I \le
\frac{\lambda(I) + \alpha(I) \lambda(I) T + (e_0(I)-(\alpha(I) +1) \lambda(I) -1)T^2
+ T^{\delta(I)+2}}{(1-T)^2}
$$

\item[(ii)]  If    $d\ge 2$ then it holds
$$
P_I \le \frac{\lambda(I) + \alpha(I) \lambda(I) T + (e_0-(\alpha(I) +1) \lambda(I) -1)T^2
+ T^{\delta(I)+2}}{(1-T)^d}.
$$

\noindent
If the equality holds then $depth(gr_I(R)) \ge d-1$.
\end{enumerate}
\end{proposition}
\begin{proof}
We write $e_i=e_i(I)$, $\lambda=\lambda(I)$,  $\alpha=\alpha(I)$ and $\delta=\delta(I)$.

\noindent
$(i)$ First we compute the right hand side power series  of  the
inequality $(i)$

\begin{multline*}
\frac{\lambda + \alpha \lambda T + (e_0-(\alpha +1) \lambda -1)T^2
+ T^{\delta+2}}{(1-T)^2}=\\
=\lambda +
\sum_{j=2}^{\delta+1}((j-1)(e_0-1) + (\alpha+2)\lambda ) T^j+
\sum_{\delta+2}^{\infty}(e_0(j+1)-e_1) T^j.
\end{multline*}
From the proof of \propref{RVV} $(i)$ we only need to prove that:

\medskip
\noindent
$(a)$ $(j-1)(e_0-1) + (\alpha+2)\lambda \ge (j-1)(e_0-1) + b +\lambda$ for all
$j=1,\cdots, \beta -1$, and

\medskip
\noindent
$(b)$ $(j-1)(e_0-1) + (\alpha+2)\lambda \ge  e_0 (j+1) - e_1$ for all
$j=\beta ,\cdots, \delta  +1$.

\medskip
\noindent
The inequality in $(a)$ is equivalent to
$(\alpha+1)\lambda \ge b $, and this inequality follows
natural epimorphism of $R/I$-modules
$(R/I)^{\mu(I)}\rightarrow I/I^2$ and the fact $\alpha(I) \ge \mu(I)-1$.

\noindent
The inequality in $(b)$ is equivalent to $\delta +1 \ge j$.

\medskip
\noindent
$(ii)$ Let $x_1,\cdots, x_{d-1} \in I$ be a superficial sequence of
$I$.
If we denote by $K=I/(x_1,\cdots,
x_{d-1})$ and $R^*=R/(x_1,\cdots, x_{d-1})$ we have
$\lambda= \length(R/I)=\length(R^{*}/K)$,
$e_i(K)=e_i(I)$, $i=0,1$, and $\delta(K)=
\delta$.

From Singh's inequality, \cite{Sin74},  we get
 $$
 P_I \le
\frac{P_{K}}{(1-T)^{d-1}},
 $$
\noindent
from $(i)$ we get the inequality in $(ii)$.
If we have
$P_I=P_{K}/(1-T)^{d-1}$,
then $depth(gr_I(R)) \ge d-1$, \cite{Sin74}.
\end{proof}


\providecommand{\bysame}{\leavevmode\hbox to3em{\hrulefill}\thinspace}


\begin{thebibliography}{10}

\bibitem{Eli90}
J.~Elias, \emph{Characterization of the {H}ilbert-{S}amuel polynomials of curve
  singularities}, Compositio Math. \textbf{74} (1990), 135--155.

\bibitem{Eli01}
\bysame, \emph{On the deep structure of the blowing-up of curve singularities},
  Math. Proc. Cambridge Philos. Soc. \textbf{131} (2001), 227--240.

\bibitem{Eli05}
\bysame, \emph{On the first normalized {H}ilbert coefficient}, J. of Pure and
  Applied Algebra \textbf{201} (2005), 116--125.

\bibitem{Lip71}
J.~Lipman, \emph{Stable ideals and {A}rf rings}, Amer. J. of Math. \textbf{93}
  (1971), 649--685.

\bibitem{Nor60}
D.G. Northcott, \emph{{A} note on the coefficients of the abstract {H}ilbert
  function}, J. London Math. Soc. \textbf{35} (1960), 209--214.

\bibitem{Ros00b}
M.~E. Rossi, \emph{A bound on the reduction number of a primary ideal}, Proc.
  A.M.S. \textbf{128} (2000), 1325--1332.

\bibitem{RV05}
M.E. Rossi and G.~Valla, \emph{The {H}ilbert function of the {R}atliff-{R}ush
  filtration}, J. of Pure and Applied Algebra \textbf{201} (2005), 25--41.

\bibitem{RVV01}
M.E. Rossi, G.~Valla, and W.V. Vasconcelos, \emph{{Maximal} {H}ilbert
  {Functions}}, Result. Math. \textbf{39} (2001), 99--114.

\bibitem{SV74}
J.~Sally and W.V. Vasconcelos, \emph{Stable rings}, J. Pure and Appl. Alg.
  \textbf{4} (1974), 319--336.

\bibitem{Sin74}
B.~Singh, \emph{Effect of a permisible blowing-up on the local {H}ilbert
  function}, Inv. Math. \textbf{26} (1974), 201--212.

\end{thebibliography}
\end{document}